\documentclass{article}
\usepackage{amsmath}
\numberwithin{equation}{section}

\begin{document}
\title{Non-recursive expressions for even-index Bernoulli numbers: A remarkable sequence of determinants.}

\author{Renaat Van Malderen}

\maketitle
\begin{abstract}Bernoulli numbers are usually expressed in terms of their lower index numbers (recursive). This paper gives explicit formulas for Bernoulli numbers of even index. The formulas contain a remarkable sequence of determinants. The value of these determinants for variable dimension is investigated.
\end{abstract}

\begin{flushleft} {\bf Keywords}:  Number theory, Bernoulli numbers, non-recursive expressions, evaluation of determinants. \end{flushleft}

\section{Introduction.}

In a previous short paper (ref. arXiv: math.NT/0503160 v1  8 Mar
2005) the author presented an explicit formula for Bernoulli numbers
of even index based on a sequence of determinants. This particular
formula was an accidental byproduct of a different mathematical
investigation. The proof via this route was fairly tedious and it
was decided to just publish the formula and see wether there was any
interest in it at all. \\

Mr Robin Chapman at the University of Exeter in the UK promptly
responded and provided a short and elegant proof using geometrical
transform techniques (generating functions). Also a request was
received from the Mathematical Society in Korea for providing a
proof for their Newsletter of JMS. Further requests were received
from China. \\

As a result this new version of the earlier short paper provides
three similar explicit expressions for Bernoulli numbers. A first
one is derived based on the sequence of determinants mentioned
earlier. The proof of the second and third one starts from the first
one and uses transform techniques similar to Mr. R.Chapman's.
Finally the value of the sequence of determinants is investigated.

\section{Non-recursive expression for even-index Bernoulli numbers.}

\begin{flushleft} First a sequence of determinants $D_{k}$ (with k= 1, 2, 3, ...) is defined as follows: \end{flushleft}

\begin{flushleft} $D_{1}$ =  $\frac{1}{3!} $ \end{flushleft}

\begin{flushleft} $D_{2}$ =  $\left|\begin{array}{cc} {\frac{1}{3!} } & {\frac{1}{5!}
} \\ {1} & {\frac{1}{3!} } \end{array}\right|$ \end{flushleft}

\begin{flushleft} $D_{3}$ = $\left|\begin{array}{ccc} {\frac{1}{3!} } & {\frac{1}{5!}
} & {\frac{1}{7!} } \\ {1} & {\frac{1}{3!} } & {\frac{1}{5!} } \\
{0} & {1} & {\frac{1}{3!} } \end{array}\right|$  \end{flushleft}

\begin{flushleft} $D_{4}$ =  $\left|\begin{array}{cccc} {\frac{1}{3!} } & {\frac{1}{5!} } & {\frac{1}{7!} } & {\frac{1}{9!} } \\ {1} & {\frac{1}{3!} } & {\frac{1}{5!} } & {\frac{1}{7!} } \\ {0} & {1} & {\frac{1}{3!} } & {\frac{1}{5!} } \\ {0} & {0} & {1} & {\frac{1}{3!} } \end{array}\right|$  etc......  \end{flushleft}

\begin{flushleft} For convenience we also define \begin{equation} \label{2.1} D_{0} = 1 \end{equation} \end{flushleft}

\begin{flushleft} Let $B_{2p}$ (with p= 1, 2, 3, ...) be the even Bernoulli number of index 2p. \end{flushleft}

\begin{flushleft} $B_{2p}$ is then given by either of the following formulas: \end{flushleft}

\begin{equation} \label{2.2} B_{2p} =\frac{(-1)^{p} (2p)!}{2} \left\{\sum _{l=0}^{p}\frac{(-1)^{l} D_{p-l} }{(2l)!} +D_{p}  \right\} \end{equation}

\begin{equation} \label{2.3} B_{2p} =-2p+\frac{(-1)^{p} (2p)!}{2^{2p} } \left\{\sum _{l=0}^{p}\frac{(-1)^{l} 3^{2l} D_{p-l} }{(2l)!} _{}
\right\} \end{equation}

\begin{equation} \label{2.4} B_{2p} =\frac{(-1)^{p+1} (2p)!D_{p} }{2(2^{2p-1} -1)} \end{equation}

\section{Relation between the determinants $D_{k}$ .}

The determinants $D_{k}$ in (2.1) are in explicit form. Nevertheless
a recursive relation among the D$_{k}$'s may be given as is obvious
from their definition. This relation may be written successively as:

\begin{flushleft} $D_{1}$ =   $\frac{1}{3!} $  $D_{0}$ \end{flushleft}

\begin{flushleft} $D_{2}$ =   $\frac{1}{3!} $  $D_{1}$ -  $\frac{1}{5!} $  $D_{0}$ \end{flushleft}

\begin{flushleft} $D_{3}$ =   $\frac{1}{3!} $  $D_{2}$ -  $\frac{1}{5!} $  $D_{1}$ +  $\frac{1}{7!} $  $D_{0}$ \end{flushleft}

\begin{flushleft} $D_{4}$ =   $\frac{1}{3!} $  $D_{3}$ -  $\frac{1}{5!} $  $D_{2}$ +  $\frac{1}{7!} $  $D_{1}$ -   $\frac{1}{9!} $  $D_{0}$ \end{flushleft}

\begin{flushleft} or in general:  \end{flushleft}

\begin{equation} \label{3.1} D_{k} =\sum _{l=1}^{k}\frac{(-1)^{l-1} D_{k-l} }{(2l+1)!} \end{equation}

\begin{flushleft} Equivalently:  \end{flushleft}

\begin{equation} \label{3.2} \sum _{l=0}^{k}\frac{(-1)^{l} D_{k-l} }{(2l+1)!}  =0  \end{equation}

\section{Proof of formula \eqref{2.2}.}

\begin{flushleft} We start from the well-known recursive relation among even-index Bernoulli numbers: \end{flushleft}

\begin{equation} \label{4.1} \sum _{k=1}^{p}\left(\begin{array}{l} {2p} \\ {2k-1} \end{array}\right)\frac{B_{2k} }{2k}  =\frac{1}{2} -\frac{1}{2p+1} \end{equation}

\begin{flushleft} The easiest way to proceed is to take a specific case, e.g. start with p=1, 2, 3: \end{flushleft}

\begin{flushleft} $$\left(\begin{array}{l} {2} \\ {1} \end{array}\right)\frac{B_{2} }{2} =\frac{1}{2} -\frac{1}{3} $$ \end{flushleft}

\begin{flushleft} $$\left(\begin{array}{l} {4} \\ {1} \end{array}\right)\frac{B_{2} }{2} +\left(\begin{array}{l} {4} \\ {3} \end{array}\right)\frac{B_{4} }{4} =\frac{1}{2} -\frac{1}{5} $$ \end{flushleft}

\begin{flushleft} $$\left(\begin{array}{l} {6} \\ {1} \end{array}\right)\frac{B_{2} }{2} +\left(\begin{array}{l} {6} \\ {3} \end{array}\right)\frac{B_{4} }{4} +\left(\begin{array}{l} {6} \\ {5} \end{array}\right)\frac{B_{6} }{6} =\frac{1}{2} -\frac{1}{7} $$ \end{flushleft}

\begin{flushleft} After some manipulation and in matrix form: \end{flushleft}

\begin{flushleft} $$\left[\begin{array}{ccc} {\frac{1}{1!} } & {0} & {0} \\ {\frac{1}{3!} } & {\frac{1}{1!} } & {0} \\ {\frac{1}{5!} } & {\frac{1}{3!} } & {\frac{1}{1!} } \end{array}\right]\left[\begin{array}{c} {\frac{B_{2} }{2!} } \\ {\frac{B_{4} }{4!} } \\ {\frac{B_{6} }{6!} } \end{array}\right]=\frac{1}{2} \left[\begin{array}{c} {\frac{1}{2!} } \\ {\frac{1}{4!} } \\ {\frac{1}{6!} } \end{array}\right]-\left[\begin{array}{c} {\frac{1}{3!} } \\ {\frac{1}{5!} } \\ {\frac{1}{7!} } \end{array}\right]$$ \end{flushleft}

\begin{flushleft} The determinant of this system equals 1.  \end{flushleft}
\pagebreak
\begin{flushleft} Solving for  $\frac{B_{6} }{6!} $  and keeping in mind the definition of the determinants $D_{k}$: \end{flushleft}

\begin{flushleft} $$\begin{array}{l} {\frac{B_{6} }{6!} =(-1)^{2} \left[\frac{1}{2} \left(\frac{D_{2} }{2!} -\frac{D_{1} }{4!} +\frac{D_{0} }{6!} \right)-\left(\frac{D_{2} }{3!} -\frac{D_{1} }{5!} +\frac{D_{0} }{7!} \right)\right]} \\ \\ {=(-1)^{3} \frac{1}{2} \left[\left(D_{3} -\frac{D_{2} }{2!} +\frac{D_{1} }{4!} -\frac{D_{0} }{6!} \right)+D_{3} \right]} \end{array}$$ \end{flushleft}

\begin{flushleft} or: \end{flushleft}

\begin{flushleft} $$B_{6} =(-1)^{3} \frac{6!}{2} \left[\sum _{l=0}^{3}\frac{(-1)^{l} D_{3-l} }{(2l)!} +D_{3}  \right]$$ \end{flushleft}

\begin{flushleft} It is obvious that this procedure may be continued for ever increasing values of p. This proves \eqref{2.2}. \end{flushleft}

\section{Proof of formulas \eqref{2.3} and \eqref{2.4}.}

Notice that formulas \eqref{2.2}, \eqref{2.3} and \eqref{2.4} also
apply to the case p=0, i.e. $B_{0}$=1. In terms of geometric
transforms, parts of formulas \eqref{2.2} and \eqref{2.3} come close
to a convolution. \\

By expressing the transforms of both  $\frac{(-1)^{l} }{(2l)!} $ and
$D_{p}$ as series in even powers of the complex variable z we
achieve convolutions. \\

Let:

\begin{flushleft} $$ D(z)=\sum _{k=0}^{\infty }D_{k} z^{2k} $$ \end{flushleft}

\begin{equation} \label{5.1} F(z)=\sum _{l=0}^{\infty }\frac{(-1)^{l} z^{2l} }{(2l+1)!}  =\frac{\sin z}{z} \end{equation}

\begin{flushleft} Convolving $D(z)$ with $F(z)$: \end{flushleft}

\begin{flushleft} $$ G(z)=D(z)F(z)=\sum _{k=0}^{\infty }\sum _{l=0}^{k}\frac{(-1)^{l} D_{k-l} z^{2k} }{(2l+1)!} $$ \end{flushleft}

\begin{flushleft} Since (see \eqref{3.2}): \end{flushleft}

\begin{flushleft}  $\sum _{l=0}^{k}\frac{(-1)^{l} D_{k-l} }{(2l+1)!}  =0$ for k=1,2... and 1 for k=0, the result is the unit impulse i.e. G(z)=1. \end{flushleft}
\pagebreak
\begin{flushleft} Therefore: \end{flushleft}

\begin{equation} \label{5.2} D(z)=\frac{z}{\sin z} \end{equation}

\begin{flushleft} Looking now at formula \eqref{2.2}, the geometric transform of \end{flushleft}

\begin{flushleft} $\displaystyle \frac{(-1)^{p} B_{2p} }{(2p)!} =\frac{1}{2} \left(\sum _{l=0}^{p}\frac{(-1)^{l} D_{p-l} }{(2l)!} +D_{p}  \right)$   equals : \end{flushleft}

\begin{equation} \label{5.3} \frac{D(z)\cos z+D(z)}{2} =D(z)\cos ^{2} \left(\frac{z}{2} \right)=\frac{z\cos \left(\frac{z}{2} \right)}{2\sin \left(\frac{z}{2} \right)} \end{equation}

\begin{flushleft} By \eqref{2.3} the transform of the same expression \end{flushleft}

\begin{flushleft} $\displaystyle \frac{(-1)^{p} B_{2p} }{(2p)!} =\frac{2p(-1)^{p+1} }{(2p)!} +\frac{1}{2^{2p} } \sum _{l=0}^{p}\frac{(-1)^{l} 3^{2l} D_{p-l} }{(2l)!}  $  equals: \end{flushleft}

\begin{flushleft} $$z\sin z+\cos \left(\frac{3z}{2} \right)D(\frac{z}{2} )=z\left(\sin z+\frac{\cos \left(\frac{3z}{2} \right)}{2\sin \left(\frac{z}{2} \right)} \right)$$ \end{flushleft}

\begin{equation} \label{5.4} =\frac{z\left(4\sin ^{2} \left(\frac{z}{2} \right)\cos \left(\frac{z}{2} \right)+4\cos ^{3} \left(\frac{z}{2} \right)-3\cos \left(\frac{z}{2} \right)\right)}{2\sin \left(\frac{z}{2} \right)} =\frac{z\cos \left(\frac{z}{2} \right)}{2\sin \left(\frac{z}{2} \right)} \end{equation}

\begin{flushleft} i.e. \eqref{5.3} and \eqref{5.4} are identical. This proves formula \eqref{2.3}. \end{flushleft}

\begin{flushleft} We proceed in a similar way for formula \eqref{2.4}: \end{flushleft}

\begin{flushleft} Let:     $W(p)=\frac{(-1)^{p} B_{2p} }{(2p)!} $ \end{flushleft}
\begin{flushleft} and:     $W(z)=\sum _{p=0}^{\infty }W_{p} z^{2p}  $  \end{flushleft}

\begin{flushleft} We already know:  \begin{equation} \label{5.5}W(z)=\frac{z\cos \left(\frac{z}{2} \right)}{2\sin \left(\frac{z}{2} \right)}
\end{equation} \end{flushleft}

\begin{flushleft} Formula \eqref{2.4} claims:  $D_{p} =-\left(2^{2p} -2\right)W_{p} $  \end{flushleft}

\begin{flushleft} The geometric transform of \eqref{5.5} equals: \end{flushleft}

\begin{flushleft} $$\frac{1}{2} \left(\frac{-2z\cos z}{\sin z} +\frac{2z\cos \left(\frac{z}{2} \right)}{\sin \left(\frac{z}{2} \right)} \right)=\frac{z}{\sin z} =D(z)$$ \end{flushleft}

\begin{flushleft} This proves \eqref{2.4}. \end{flushleft}

\section{Value of the determinants $D_{k}$.}

From \eqref{5.2}:

\begin{equation} \label{6.1} D(z)=\frac{z}{\sin z} =\sum_{p=0}^{\infty }D_{p} z^{2p} \end{equation}

\begin{equation} \label{6.2} D_{p} =\frac{1}{2\pi i} \oint \frac{D(z)dz}{z^{2p+1} }  =\frac{1}{2\pi i} \oint \frac{dz}{z^{2p} \sin z}  \end{equation}

\begin{flushleft} where the integral is taken along a sufficiently small circle around the origin z=0. \end{flushleft}

\begin{flushleft} We now use a well-known result: \end{flushleft}

\begin{flushleft} $$\frac{1}{\sin z} =\frac{1}{z} +2z\sum _{k=1}^{\infty }\frac{(-1)^{k+1} }{\pi ^{2} k^{2} } \cdot \frac{1}{\left(1-\frac{z^{2} }{\pi ^{2} k^{2} } \right)} $$ \end{flushleft}

\begin{flushleft} For  $\left|z\right|\langle \pi $  and k=1,2,.... , we are allowed to expand the last term in a power series. By dividing further by z$^{2p}$ we obtain: \end{flushleft}

\begin{equation} \label{6.3} \frac{1}{z^{2p} \sin z} =\frac{1}{z^{2p+1} } +2z\sum _{k=1}^{\infty }\frac{(-1)^{k+1} }{\pi ^{2} k^{2} } \cdot \sum _{l=0}^{\infty }\frac{1}{\pi ^{2l} k^{2l} z^{2(p-l)-1} } \end{equation}

\begin{flushleft} The value of \eqref{6.2} equals the residue of \eqref{6.3}. \end{flushleft}

\begin{flushleft} For p=0 the residue equals 1, i.e. $D_{0}$ = 1 which is correct. \end{flushleft}

\begin{flushleft} For p $>=$ 1 we look for terms where 2p-2l-1 = 1 or l=p-1, resulting in: \end{flushleft}

\begin{equation}\label{6.4} D_{p} =\frac{2}{\pi ^{2p} } \sum _{k=1}^{\infty }\frac{(-1)^{k+1} }{k^{2p} } =\frac{2}{\pi ^{2p} } \left(1-\frac{1}{2^{2p-1} } \right)\zeta (2p) \end{equation}

\begin{flushleft} or \end{flushleft}

\begin{equation}\label{6.5} \zeta (2p)=\frac{\pi ^{2p} D_{p} }{2(1-2^{1-2p} )} \end{equation}

\begin{flushleft} in case one prefers to express it the other way around. \end{flushleft}
\pagebreak
\begin{flushleft} It is easily seen from \eqref{6.4} that for large p: \end{flushleft}

\begin{equation}\label{6.6} D_{p} \cong \frac{2}{\pi ^{2p} }  \end{equation}

\begin{flushleft} or: \end{flushleft}

\begin{flushleft} $$\mathop{\lim }\limits_{p\to \infty } =\frac{\pi ^{2p} D_{p} }{2} = 1$$ \end{flushleft}

\begin{flushleft}  \end{flushleft}
\begin{flushleft}  \end{flushleft}
\begin{flushleft}
Renaat Van Malderen \\ Address: Maxlaan 21, 2640 Mortsel, Belgium \\
Author can be reached by email: hans.vanmalderen@bertholdtech.com
\end{flushleft}

\end{document}